%% file: rado.tex
\documentclass[12pt]{amsart}
\usepackage{epsfig,color}
\usepackage{amsfonts, amsthm, amsmath}

\theoremstyle{plain}
\newtheorem{theorem}{Theorem}

\newtheorem{lemma}[theorem]{Lemma}
\newtheorem{corollary}[theorem]{Corollary}
\newtheorem{proposition}[theorem]{Proposition}

\begin{document}

\title[A generalization of Rado for almost graphical boundaries]{A generalization of Rado's Theorem for almost graphical boundaries}
\author{Brian Dean
    \and
    Giuseppe Tinaglia}
\address{Department of Mathematics\\
    Hylan Building\\
    University of Rochester\\
    Rochester, NY  14627}
\email{bdean@math.rochester.edu}
\address{Department of Mathematics\\
    Johns Hopkins University\\
    3400 N. Charles St.\\
    Baltimore, MD  21218}
\email{tinaglia@math.jhu.edu}
\date{}

\begin{abstract}
In this paper, we prove a generalization of Rado's Theorem, a
fundamental result of minimal surface theory, which says that
minimal surfaces over a convex domain with graphical boundaries
must be disks which are themselves graphical.  We will show that,
for a minimal surface of any genus, whose boundary is ``almost
graphical'' in some sense, that the surface must be graphical once
we move sufficiently far from the boundary.
\end{abstract}

\maketitle

\section{Introduction}\label{intro}

One of the fundamental results of minimal surface theory is Rado's
Theorem, which is connected to the famous Plateau Problem.  Rado's
Theorem (see \cite{Ra}) states that if $\Omega\subset\mathbb{R}^2$
is a convex subset and $\sigma\subset\mathbb{R}^3$ is a simple
closed curve which is graphical over $\partial\Omega$, then any
minimal surface $\Sigma\subset\mathbb{R}^3$ with
$\partial\Sigma=\sigma$ must be a disk which is graphical over
$\Omega$, and hence unique by the maximum principle.

The proof of Rado's Theorem begins by assuming there is a point at
which $\Sigma$ is not graphical, i.e., where $\Sigma$ has a
vertical tangent plane.  One can then use the description of the
local intersection of minimal surfaces as $n$-prong singularities
(see, for example, \cite[Section 4.6]{CM2}) to derive a
contradiction of the assumption on the boundary $\sigma$.

Our goal is to generalize Rado's Theorem for the case in which
$\sigma$ is not graphical, but satisfies some ``almost graphical''
condition.  We will show that, although $\Sigma$ is obviously not graphical
near its boundary, if we move far enough in from the boundary,
$\Sigma$ will be graphical.  We will prove the genus zero case
first, and then generalize to higher genus surfaces.

Throughout this paper, we will use the topological fact that, if
$\Sigma$ has genus $n$, any collection of $n+1$ disjoint closed
simple curves on $\Sigma$ must separate $\Sigma$ into at least two
connected components.  This can be seen as follows.  If a compact
connected orientable surface with boundary of genus $n$ has $k$
boundary components, its Euler characteristic is $2-2n-k$.  By
cutting the surface along a circle which does not separate the
surface into at least two connected components, the number of
boundary components would increase by two, while the Euler
characteristic would stay the same.  Thus, the genus would
decrease by one.  Therefore, the maximum number of such cuts would
be equal to the genus $n$, so any collection of $n+1$ disjoint
closed simple curves must separate the surface into at least two
connected components.

The ``almost graphical'' condition we will use will be as follows.
Let $\sigma=(\sigma_1,\sigma_2,\sigma_3)$ be a parametrization for
$\partial\Sigma$.  Then, we say that $\sigma$ is ``$C,h$-almost
graphical'' if

1) $\sigma$ has one connected component.

2) After possibly a rotation, $|\sigma_3|<Ch$ (and thus, all of
$\Sigma$ lies in a narrow vertical slab).

3) $\sigma$ is ``$h$-almost monotone'', i.e., for any
$y\in\sigma$, $B_{4h}(y)\cap\sigma$ has only one component which
intersects $B_{2h}(y)$.  Therefore, any point $b\in
B_{2h}(y)\cap\sigma$ can be joined to $y$ by a path in
$B_{4h}(y)\cap\sigma$.  See Figure~\ref{figure00}.

\begin{figure}[h]
\vskip 15pt
    \setlength{\captionindent}{4pt}
    \centering\input{figu00.pstex_t}
   \caption{}\label{figure00}\vskip 15pt
\end{figure}

We now state our main result.

\begin{theorem}\label{mainresult}
There exists a $C>0$ (not depending on $\Sigma$) such that if $\Sigma$ is an embedded minimal
surface of genus $n$, $n\ge 0$, with $C,h$-almost graphical boundary $\sigma
=\partial\Sigma\subset \partial B_R$, then $\Sigma \cap B_{R-(64n+30)h}$
is graphical.
\end{theorem}

\section{Catenoid foliations}\label{catenoid}

Although the proof of Rado's Theorem utilizes intersections of
minimal surfaces with planes (namely, vertical tangent planes),
the proof of our generalization will require a greater degree of
sophistication.  We will be intersecting minimal surfaces at
nongraphical points with carefully chosen catenoids.

Here, we provide some background and important results involving
catenoid foliations.  This material is covered in greater detail
(including proofs) in \cite[Appendix A]{CM1}.  Recall that in this
paper we will be talking about minimal surfaces $\Sigma$ which lie
in a narrow vertical slab.

\begin{figure}[h]
\vskip 15pt
    \setlength{\captionindent}{4pt}
    \centering\input{fig000.pstex_t}
   \caption{The catenoid foliation}\label{fig000}\vskip 15pt
\end{figure}
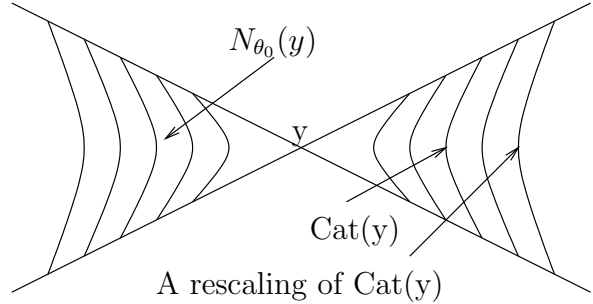

Let $\mbox{Cat}(y)$ be the vertical catenoid centered at
$y=(y_1,y_2,y_3)\in\mathbb{R}^3$.  In other words,
$$\mbox{Cat}(y)=\{x\in\mathbb{R}^3|\cosh^2(x_3-y_3)=(x_1-y_1)^2+(x_2-y_2)^2\}.$$
For an angle $\theta\in \left(0,\frac{\pi}{2}\right)$, we denote
by $\partial N_{\theta}(y)$ the cone
$$\{x\in\mathbb{R}^3|(x_3-y_3)^2=|x-y|^2\sin^2\theta\}.$$
Then, we see that $\partial N_{\pi/4}(y)\cap
\mbox{Cat}(y)=\emptyset$, since $\cosh t>t$ for all $t\ge 0$.  So,
we set
$$\theta_0=\inf\{\theta|\partial N_{\theta}(y)\cap\mbox{Cat}(y)=\emptyset\}.$$
Thus, $\partial N_{\theta_0}(y)$ and $\mbox{Cat}(y)$ intersect
tangentially in a pair of circles (one above the $x_1x_2$-plane
and one below).  Let $\mbox{Cat}_0(y)$ be the component of
$\mbox{Cat}(y)\backslash\partial N_{\theta_0}(y)$ containing the
neck
$$\{x|x_3=y_3, (x_1-y_1)^2+(x_2-y_2)^2=1\}.$$
If $x\in\mbox{Cat}_0(y)$, then the line segment joining $y$ and
$x$ intersects $\mbox{Cat}_0(y)$ at exactly one point, namely $x$.
So, the dilations of $\mbox{Cat}_0(y)$ about $y$ are disjoint, and
give us a minimal foliation (see Figure \ref{fig000})
of the solid (open) cone
$$N_{\theta_0}(y)=\{x\in\mathbb{R}^3|(x_3-y_3)^2<|x-y|^2\sin^2\theta_0\}.$$
The leaves of this foliation all have boundary in $\partial
N_{\theta_0}(y)$ and are the level sets of the function $f_y$
given by
$$y+\frac{x-y}{f_y(x)}\in\mbox{Cat}_0(y).$$
Choose $\beta_A>0$ sufficiently small so that
$$\{x||x_3-y_3|\le 2\beta_A h\}\backslash B_{h/8}(y)\subset N_{\theta_0}(y)$$
and
\begin{equation}\label{betacondition}
\{x|f_y(x)=3h/16\}\cap\{x||x_3-y_3|\le 2\beta_Ah\}\subset
B_{7h/32}(y).
\end{equation}
Since the intersection of any two minimal surfaces (in particular,
our given $\Sigma$ and any catenoid in our foliation) is locally
given by an $n$-prong singularity, i.e., $2n$ embedded arcs which
meet at equal angles (see Claim 1 of Lemma 4 in \cite{HoMe}), we
get the following Lemma.

\begin{lemma}\label{catenoidnpronglemma}
If $z\in\Sigma\subset N_{\theta_0(y)}$ is a nontrivial interior
critical point of $f_y|_{\Sigma}$, then $\{x\in\Sigma
|f_y(x)=f_y(z)\}$ has an $n$-prong singularity at $z$ with $n\ge
2$.
\end{lemma}

As a consequence of this, we obtain a version of the Strong
Maximum Principle.

\begin{lemma}\label{catenoidmaxprinciplelemma}
If $\Sigma\subset N_{\theta_0}(y)$, then $f_y|_{\Sigma}$ has no
nontrivial interior local extrema.
\end{lemma}

Using this, we can show, using the foliation indexed by $f_y$,
that a minimal surface in a narrow slab either stays near the
boundary or comes near the center.

\begin{corollary}\label{corollarytocatenoidmaxprinciplelemma}
If $\partial\Sigma\subset\partial B_h(y)$,
$B_{3h/4}(y)\cap\Sigma\ne\emptyset$, and
$$\Sigma\subset B_h(y)\cap\{x||x_3-y_3|\le 2\beta_A h\},$$
then $B_{h/4}(y)\cap\Sigma\ne\emptyset$.
\end{corollary}

Iterating Corollary~\ref{corollarytocatenoidmaxprinciplelemma}
along a chain of balls, we see that we will be able to extend
curves out close to the boundary.  Here, $T_h$ refers to a tubular
neighborhood of radius $h$, and $\gamma_{p,q}$ is the line segment
joining $p$ and $q$.

\begin{corollary}\label{curveextensioncorollary}
If $\Sigma\subset\{|x_3|\le 2\beta_A h\}$, points
$p,q\in\{x_3=0\}$ satisfy
$T_h(\gamma_{p,q})\cap\partial\Sigma=\emptyset$, and
$$y_p\in B_{h/4}(p)\cap\Sigma,$$
then a curve $\nu\subset T_h(\gamma_{p,q})\cap\Sigma$ connects
$y_p$ to $B_{h/4}(q)\cap\Sigma$.
\end{corollary}

The final catenoid foliation result we will use shows that the
vertical projection of $\Sigma$ cannot stray too far outside the
vertical projection of $\partial\Sigma$.

\begin{corollary}\label{unboundedcomponentprojectionlemma}
If $\Sigma\subset\{|x_3|\le 2\beta_A h\}$ and $E$ is an unbounded
component of
$$\mathbb{R}^2\backslash T_{h/4}(\Pi(\partial\Sigma)),$$
then $\Pi(\Sigma)\cap E=\emptyset$.
\end{corollary}

\section{The genus zero case}\label{genuszero}

\noindent{\it Proof of Theorem~\ref{mainresult} for $n=0$}.  Let
$C=\beta_A$, where $\beta_A>0$ is defined by
(\ref{betacondition}).  Let $\Sigma$ be an embedded minimal disk
with $C,h$-almost graphical boundary $\sigma
=\partial\Sigma\subset\partial B_R$.  Our proof begins with an
argument which is similar to the first step of the proof of
\cite[Lemma I.0.11]{CM1}.  Suppose $\Sigma\cap B_{R-30h}$ is not
graphical; let $z\in\Sigma\cap B_{R-30h}$ be a point such that the
tangent plane to $\Sigma$ at $z$ is vertical.  Fix $y\in\partial
B_{4h}(z)$ so that the line segment $\gamma_{y,z}$ is normal to
$\Sigma$ at $z$.  Then, $f_y(z)=4h$, where $f_y$ is the function
used to define the level sets of the catenoid foliation for
catenoids centered at $y$; see Section~\ref{catenoid}.  Let $y'$
be given such that $y'\in\partial B_{10h}(y)$ and
$z\in\gamma_{y,y'}$.

Any simple closed curve $\rho\subset\Sigma\backslash\{f_y >4h\}$
bounds a disk $\Sigma_{\rho}\subset\Sigma$.  By
Lemma~\ref{catenoidmaxprinciplelemma}, $f_y$ has no maxima on
$\Sigma_{\rho}\cap \{f_y >4h\}$ so that $\Sigma_{\rho}\cap\{f_y
>4h\}=\emptyset$.  On the other hand, by
Lemma~\ref{catenoidnpronglemma}, there is a neighborhood
$U_z\subset\Sigma$ of $z$ such that $U_z\cap\{f_y =4h\}$ is an
$n$-prong singularity with $n\ge 2$; in other words, $U_z\cap\{f_y
=4h\}\backslash \{z\}$ is the union of $2n\ge 4$ disjoint embedded
arcs meeting at $z$.  Moreover, $U_z\backslash\{f_y \ge 4h\}$
(i.e., the part of $U_z$ inside the catenoid $\{f_y =4h\}$) has
$n$ components $U_1,\ldots,U_n$ with
$$\overline{U_i}\cap\overline{U_j}=\{z\}\;\mbox{for}\; i\ne j.$$
If a simple curve $\widetilde{\rho_z}\subset\Sigma\backslash\{f_y
\ge 4h\}$ connects $U_1$ to $U_2$, connecting
$\partial\widetilde{\rho_z}$ by a curve in $U_z$ gives a simple
closed curve $\rho_z\subset\Sigma\backslash \{f_y >4h\}$ with
$\widetilde{\rho_z}\subset \rho_z$ and $\rho_z\cap\{f_y \ge
4h\}=\{z\}$.  Hence, $\rho_z$ bounds a disk
$\Sigma_{\rho_z}\subset\Sigma\backslash\{f_y >4h\}$.  By
construction,
$$U_z\cap\Sigma_{\rho_z}\backslash\cup_i \overline{U_i}\ne\emptyset.$$
This is a contradiction, so $U_1$ and $U_2$ must be contained in
components $\Sigma_{4h}^1\ne\Sigma_{4h}^2$ of
$\Sigma\backslash\{f_y \ge 4h\}$ with
$z\in\overline{\Sigma_{4h}^1}\cap\overline{\Sigma_{4h}^2}$.  For
$i=1,2$, Lemma~\ref{catenoidmaxprinciplelemma} and
(\ref{betacondition}) give us $y_i^a\in
B_{h/4}(y)\cap\Sigma_{4h}^i$.  By
Corollary~\ref{curveextensioncorollary}, there exist curves
$\nu_i\subset T_h(\gamma_{y,y'})\cap\Sigma$ with
$\partial\nu_i=\{y_i^a,y_i^b\}$, where $y_i^b\in B_{h/4}(y')$.
There are two cases:
\begin{itemize}
\item If $y_1^b$ and $y_2^b$ do not connect in
$B_{4h}(y')\cap\Sigma$, take $\gamma_0\subset B_{5h}(y)\cap\Sigma$
from $y_1^a$ to $y_2^a$ and set
$\gamma_a=\nu_1\cup\gamma_0\cup\nu_2$ and $y_i=y_i^b$. \item
Otherwise, if $\widehat{\gamma}_0\subset B_{4h}(y')\cap\Sigma$
connects $y_1^b$ and $y_2^b$, set
$\gamma_a=\nu_1\cup\widehat{\gamma}_0\cup\nu_2$ and $y_i=y_i^a$.
\end{itemize}
In either case, after possibly switching $y$ and $y'$, we get a
curve
$$\gamma_a\subset (T_h(\gamma_{y,y'})\cup B_{5h}(y'))\cap\Sigma$$
with $\partial\gamma_a=\{y_1,y_2\}\subset B_{h/4}(y)$ and $y_i\in
S_i^a$ for components $S_1^a\ne S_2^a$ of $B_{4h}(y)\cap\Sigma$.

Let $y''$ be given so that $y''\in\partial B_{\sqrt{R^2-h^2}}$ and
$y\in\gamma_{y',y''}$ (that is, $y'$, $z$, $y$, and $y''$ are all
collinear).  By Corollary~\ref{curveextensioncorollary}, for
$i=1,2$, there exist curves $\mu_i\subset
T_h(\gamma_{y,y''})\cap\Sigma$ connecting $y_i$ to points $z_i\in
B_{h/4}(y'')\cap\Sigma$.  We then intersect $\Sigma$ with a plane
connecting the points $z_1$ and $z_2$, and which is perpendicular
to $\gamma_{y,y''}$ %
%
%
%
%
(such a plane exists since $z_i$ can be taken to be any point in
$B_{h/4}(y'')\cap\mu_i$). If $z_1$ and $z_2$ can be connected by a
curve $\lambda$ in the intersection of $\Sigma$ with this plane,
then $y_1$ and $y_2$ are connected by the curve
$\mu_1\cup\lambda\cup\mu_2$.  This contradicts the following
Claim.

\bigskip
\noindent\textbf{Claim:} Let $H$ be the half-space given by
$$H=\{x|\langle y-y',x-y\rangle>0\}.$$
Then, $y_1$ and $y_2$ can not be connected by a curve in
$T_h(H)\cap\Sigma$.

\bigskip
\noindent\textit{Proof of Claim} \cite[p. 11-12]{CM1}.  If
$\eta_{1,2}\subset T_h(H)\cap\Sigma$ connects $y_1$ and $y_2$,
then $\eta_{1,2}\cup\gamma_a$ bounds a disk
$\Sigma_{1,2}\subset\Sigma$.  Since $\eta_{1,2}\subset T_h(H)$,
$\partial B_{8h}(y')\cap\partial\Sigma_{1,2}$ consists of an odd
number of points in each $S_i^a$ and hence $\partial
B_{8h}(y')\cap\Sigma_{1,2}$ contains a curve from $S_1^a$ to
$S_2^a$.  However, $S_1^a$ and $S_2^a$ are distinct components of
$B_{4h}(y)\cap\Sigma$, so this curve must contain a point
$$y_{1,2}\in\partial B_{4h}(y)\cap\partial B_{8h}(y')\cap\Sigma_{1,2}.$$
By construction, $\Pi (y_{1,2})$ is in an unbounded component of
$\mathbb{R}^2\backslash T_{h/4}(\Pi (\partial\Sigma_{1,2}))$,
contradicting Corollary~\ref{unboundedcomponentprojectionlemma}.
Therefore, $y_1$ and $y_2$ can't be connected in
$T_h(H)\cap\Sigma$.\qed

\bigskip
Returning to the proof of the genus zero case of Theorem~\ref{mainresult}, since
there is no curve in the intersection of $\Sigma$ and this plane
which connects $z_1$ to $z_2$, we get disjoint curves $\lambda_1$
and $\lambda_2$ in the intersection of $\Sigma$ with this plane,
with $z_i\in\lambda_i$ for $i=1,2$.  Neither $\lambda_1$ nor
$\lambda_2$ can be closed, as if either were closed, it would
bound a disk in the intersection of $\Sigma$ and the plane,
violating the maximum principle.  So, these curves must go to the
boundary of $\Sigma$, i.e., there exist points
$b_i\in\lambda_i\cap\sigma$ for $i=1,2$.  By construction, $b_2\in
B_{2h}(b_1)$.  By the $h$-almost monotonicity of $\sigma$, there
is a curve %
%
%
%
$\alpha\subset B_{8h}(b_1)\cap\sigma$ connecting $b_1$ and $b_2$.
Thus, $y_1$ is connected to $y_2$ by the curve
$\mu_1\cup\lambda_1\cup\alpha\cup\lambda_2\cup\mu_2$,
contradicting that $y_1$ can not be connected to $y_2$ in
$T_h(H)\cap\Sigma$.  Therefore, $\Sigma\cap B_{R-30h}$ is
graphical.\qed

\section{The higher genus case}\label{highergenus}

In this section we prove the higher genus case, i.e.,
Theorem~\ref{mainresult} for $n\ge 1$.   We begin by looking at the
case $n=1$, that is
\begin{theorem}\label{mainresultgenusone}
Let $\Sigma$ be an embedded minimal surface with genus 1 such that
$\sigma =\partial\Sigma\subset \partial B_{R}$ is $C,h$-almost
graphical.  Then, $\Sigma\cap B_{R-94h}$ is graphical.
\end{theorem}

As in the genus 0 case, when we say $\partial\Sigma$ is
$C,h$-almost graphical, we will be taking $C=\beta_A$.

In dealing with the genus zero case, each time we had a closed
path we could  say that it bounded a disk; that is, each closed
path was homotopic to a point. This is no longer true when the
genus is one or higher.  However, given any genus one surface, any
two disjoint closed paths divide the surface into at least two
regions.  In the following lemma, we show that, for our given
minimal surface $\Sigma$, it is impossible to have two nontrivial
closed paths (i.e., two closed paths which are not homotopic to a
point) which are far apart.

\begin{lemma}\label{twodistantpaths}
There can't be two nontrivial closed simple paths $\gamma_1$  and
$\gamma_2$ in $\Sigma$ so that
\begin{equation}\label{twodstanpathseq}\text{dist}(T_{\frac{h}{4}}(\gamma_1),T_{\frac{h}{4}}(\gamma_2))>0.\end{equation}
\end{lemma}

\begin{proof}
Because of the topological properties of a genus one surface,
$\gamma_1\cup\gamma_2$ bounds a connected region
$\Sigma'\subset\Sigma$.  See Figure~\ref{figure0}, which shows a
closed torus (the same is true for any surface of genus one).
\begin{figure}[h]
\vskip 15pt
    \setlength{\captionindent}{4pt}
    \centering\input{figu0.pstex_t}
   \caption{}\label{figure0}\vskip 15pt
\end{figure} However, this topological fact and \eqref{twodstanpathseq} gives $\mathbb{R}^2\backslash
T_{h/4}(\Pi(\partial\Sigma'))\neq\emptyset$. This contradicts
Corollary \ref{unboundedcomponentprojectionlemma}.
\end{proof}

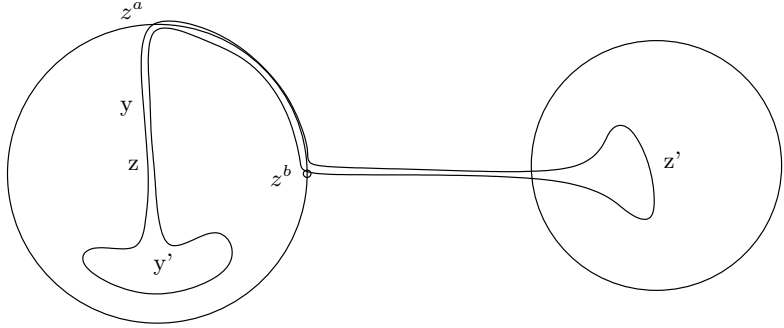
\begin{figure}[h]
\vskip 15pt
    \setlength{\captionindent}{4pt}
    \centering\input{figu1.pstex_t}
   \caption{Path in a dumbbell}\label{figure1}\vskip 15pt
\end{figure}

To apply Lemma \ref{twodistantpaths} we need to be able to build a
closed path which is not homotopic to a point and is contained in
a fixed region. In the following lemma we build a path which is
not homotopic to a point and is contained in a dumbbell-shaped
region as shown in Figure~\ref{figure1}.

\begin{lemma}\label{dumbbellshapedpath}
Let $\Sigma$ be as in Theorem \ref{mainresult} with genus $n\ge 1$, and let
$z\in B_{R-30h}\cap\Sigma$ be such that $z$ is not graphical.
Then, for any $z'\in \partial \Sigma$, there exists a closed path
contained in $(B_{8h}(z')\cup T_h(\gamma_{z,z'})\cup
B_{13h}(z))\cap \Sigma$ which is not homotopic to a point.
\end{lemma}

\begin{proof}
The construction of this path starts in the same way as the
construction at the beginning of the proof of the genus zero case
of Theorem~\ref{mainresult}. Fix $y\in\partial B_{4h}(z)$ so that
the line segment $\gamma_{y,z}$ is normal to $\Sigma$ at $z$.
Then, $f_y(z)=4h$, where $f_y$ is the function used to define the
level sets of the catenoid foliation for catenoids centered at
$y$.  Now, by Lemma~\ref{catenoidnpronglemma}, there is a
neighborhood $U_z\subset\Sigma$ of $z$ such that $U_z\cap\{f_y
=4h\}$ is an $n$-prong singularity with $n\ge 2$; in other words,
$U_z\cap\{f_y =4h\}\backslash \{z\}$ is the union of $2n\ge 4$
disjoint embedded arcs meeting at $z$.  Moreover,
$U_z\backslash\{f_y \ge 4h\}$ (i.e., the part of $U_z$ inside the
catenoid $\{f_y =4h\}$) has $n$ components $U_1,\ldots,U_n$ with
$$\overline{U_i}\cap\overline{U_j}=\{z\}\;\mbox{for}\; i\ne j.$$

However, unlike in the genus zero case, we can not say that that
$U_1$ and $U_2$ must be contained in distinct components
$\Sigma_{4h}^1\ne\Sigma_{4h}^2$ of $\Sigma\backslash\{f_y \ge
4h\}$ with
$z\in\overline{\Sigma_{4h}^1}\cap\overline{\Sigma_{4h}^2}$.  There
are two possibilities: either the components are distinct (as in
the genus zero case), or they coincide, i.e.,
$\Sigma_{4h}^1=\Sigma_{4h}^2$.  First, we consider the case where
the components coincide.  Then, as in the proof of the genus zero
case, we would have a simple curve
$\widetilde{\rho_z}\subset\Sigma\backslash\{f_y\ge 4h\}$
connecting $U_1$ to $U_2$, and connecting
$\partial\widetilde{\rho_z}$ by a curve in $U_z$, we can build a
simple closed curve $\rho_z\subset\Sigma\backslash\{f_y >4h\}$
with $\widetilde{\rho_z}\subset\rho_z$ and $\rho_z\cap\{f_y\ge
4h\}=\{z\}$.  We saw in the proof of the genus zero case 
that $\rho_z$ can not bound a
disk.  Thus, the curve $\rho_z$ is not homotopic to a point, and
so the lemma is proved in this case, since $\rho_z\subset
\Sigma\cap\{f_y\le 4h\}\subset \Sigma\cap B_{13h}(z)$.

If instead the components $\Sigma_{4h}^1$ and $\Sigma_{4h}^2$ are
distinct (as in the genus zero case), then let $y'$ be given such
that $y'\in\partial B_{10h}(y)$ and $z\in \gamma_{y,y'}$.  Then,
as in the proof of the genus zero case, we get a
curve
$$\gamma_a\subset (T_h(\gamma_{y,y'})\cup B_{5h}(y'))\cap\Sigma$$
with $\partial\gamma_a=\{y_1,y_2\}\subset B_{h/4}(y)$ and $y_i\in
S_i^a$ for components $S_1^a\ne S_2^a$ of $B_{4h}(y)\cap\Sigma$.

Let $S$ be the annulus given by $B_{11h}(z)\backslash B_{10h}(z)$,
let $z^a\in\partial B_{11h}(z)$ such that $y\in\gamma_{z,z^a}$ and
let $z^b\in\partial B_{11h}(z)$ such that $z^b\in\gamma_{z,z'}$.
Let $\bar{\gamma}$ be the shortest polygonal path in $S$
connecting $z^a$ and $z^b$ then, using
Corollary~\ref{curveextensioncorollary} we can first build two
paths $\gamma^a_i$ connecting $y_i$ to $z^a_i\in
B_\frac{h}{4}(z^a)$ with $\gamma^a_i\subset T_h(\gamma_{y,z^a})$.
Then we can build another two paths $\gamma^{ab}_i$ from $z^a_i$
to $z^b_i\in B_\frac{h}{4}(z^b)$ with $\gamma^{ab}_i\subset
T_h(\bar{\gamma})$. Eventually, we build two paths $\gamma^b_i$
connecting $z^b_i$ to $z'_i\in B_\frac{h}{4}(\bar{z}')$ with
$\gamma^{b}_i\subset T_h(\gamma_{z^b,\bar{z}'})$, where
$\bar{z}'\in\partial B_{\sqrt{R^2-h^2}}$ such that
$\bar{z}'\in\gamma_{z^b,z'}$. So far we have a path
$\gamma=\gamma^b_1\cup\gamma^{ab}_1\cup\gamma^a_1\cup\gamma_a\cup\gamma^a_2\cup\gamma^{ab}_2\cup\gamma^b_2$
whose ends, $z'_i$, are contained in $B_\frac{h}{4}(\bar{z}')$. As
in the end of the proof of the genus zero case, taking the
intersection of $\Sigma$ with a plane connecting the points $z'_1$
and $z'_2$, and which is perpendicular to $\gamma_{z',z^b}$ we can
connect $z'_1$ and $z'_2$ with a path, $\gamma_b$, that is
contained in $B_{8h}(z')\cap\partial\Sigma$. We are left to show
why
$\gamma=\gamma^b_1\cup\gamma^{ab}_1\cup\gamma^a_1\cup\gamma_a\cup\gamma^a_2\cup\gamma^{ab}_2\cup\gamma^b_2\cup\gamma_b$
is not homotopic to a point.

In the proof of the genus zero case, in the claim, we proved that
if a path $\eta_{1,2}\subset T_h(H)\cap\Sigma$ where
$H=\{x|\langle y-y', x-y \rangle >0 \}$ connects $y_1$ and $y_2$,
the ends of $\gamma_a$, then $\gamma_a\cup\eta_{1,2}$ can not
bound a disk. However, a closer look at the proof of that shows
that as long as $\partial B_{4h}(y)\cap\partial
B_{8h}(y')\cap\Sigma$ is in the unbounded component of
$\mathbb{R}^2\backslash T_{h/4}(\Pi(\gamma_a\cup\eta_{1,2}))$ then
$\gamma_a\cup\eta_{1,2}$ cannot bound a disk. Since this is the
case, $\gamma$ cannot bound a disk.
\end{proof}
As a consequence of Lemma \ref{twodistantpaths} and Lemma
\ref{dumbbellshapedpath} we have the following lemma that says
that if the interior of the surface fails to be graphical at two
points, then these two points have to be close. In other words,
the interior minus a smaller ball is graphical.

\begin{lemma}\label{graphminusaball}
Let $\Sigma$ be as in Theorem~\ref{mainresult} with genus $n\ge 1$, and let
$z_1, z_2\in B_{R-30h}\cap\Sigma$ such that $z_1$ and $z_2$ are
not graphical.  Then, $|z_1-z_2|<31h$.
\end{lemma}
\begin{proof}
Assume $|z_1-z_2|\geq31h$ and let $\pi$ be the plane perpendicular
to $\gamma_{z_1,z_2}$ through its midpoint. Fix
$z_1',z_2'\in\partial\Sigma$ so that $z_i'$ is in the half space
containing $z_i$ and $|z_1'-z_2'|>17h$. Then applying Lemma
\ref{dumbbellshapedpath} creates two paths $\gamma_1$ and
$\gamma_2$ as in Lemma~\ref{twodistantpaths}, see
Figure~\ref{figure2}, giving the contradiction.
\end{proof}

\begin{figure}[h]
\vskip 15pt
    \setlength{\captionindent}{4pt}
    \centering\input{figu2.pstex_t}
   \caption{}\label{figure2}\vskip 15pt
\end{figure}

We are now ready to prove Theorem~\ref{mainresultgenusone}.

\begin{proof}[Proof of Theorem~\ref{mainresultgenusone}.]
Lemma~\ref{graphminusaball} says that if $z_1$ and $z_2$ are two
nongraphical points in $B_{R-30h}$ then $(B_{R-30h}\backslash
B_{31h}(z_1))\cap\Sigma$ is graphical. If $z_1\in B_{R-30h-32h}$
then the annulus $(B_{R-30h}\backslash B_{R-31h})\cap \Sigma$ is
graphical and applying Rado's theorem gives that
$B_{R-30h}\cap\Sigma$ is graphical. If instead $z_1$ is not in
$B_{R-62h}$ then $B_{R-94h}\cap\Sigma$ is graphical. In either
case, the theorem follows.
\end{proof}

Now we begin to prove the general case of Theorem~\ref{mainresult}.

Let $A_i=B_{R-64ih-30h}\backslash B_{R-64ih-31h}$ for $i=0,...,n$.
These $n+1$ annuli have width $h$, and the distance between $A_i$
and $A_{i+1}$ is $64h$. Theorem~\ref{mainresult} will
clearly follow once we have proved the following proposition,

\begin{proposition}\label{graphicalannulus}
There exists an $i=0,...,n$ such that $A_i\cap\Sigma$ is
graphical.
\end{proposition}

The proof of Proposition~\ref{graphicalannulus} uses the
equivalent of Lemma~\ref{twodistantpaths} for the genus $n$ case.

\begin{lemma}\label{ndistantpaths} There can't be $n+1$ nontrivial
closed simple paths $\gamma_i\in\Sigma$, $i=0,...,n$, so that%
\begin{equation}\label{ndstanpathseq}\text{dist}(T_{\frac{h}{4}}(\gamma_i),T_{\frac{h}{4}}(\gamma_j))>0, \quad \text{for } i\neq j.\end{equation}
\end{lemma}

\begin{proof}[Proof of Lemma~\ref{ndistantpaths}]
The proof uses the same idea that it is used to prove
Lemma~\ref{twodistantpaths}. In the genus n case we use the
topological property that $n+1$ disjoint closed simple paths bound
at least one connected region with more than one boundary
component.
\end{proof}

\begin{figure}[!h] \vskip 15pt
    \setlength{\captionindent}{4pt}
    \centering\input{figu3.pstex_t}
    \caption{}\label{figure3}\vskip 15pt
\end{figure}

\begin{proof}[Proof of Proposition~\ref{graphicalannulus}]
Let us assume the proposition is false and let $z_i\in
A_i\cap\Sigma$ be nongraphical points. Working as we did in the
proof of Lemma~\ref{dumbbellshapedpath} we can find
$z_i'\in\partial\Sigma$ and polygonal paths $\bar{\gamma}_i$, and
build $n+1$ nontrivial paths $\gamma_i\subset B_{31h}(z_i)\cup
T_h(\bar{\gamma_i}) \cup B_{8h}(z_i')$ (where the $n+1$ boundary
points $z_i'$ are chosen in such a way that the corresponding
dumbbell-shaped regions are pairwise disjoint) such that %
\eqref{ndstanpathseq} in Lemma~\ref{ndistantpaths} holds.  See
Figure~\ref{figure3}.

This contradicts Lemma~\ref{ndistantpaths}, proves
Proposition~\ref{graphicalannulus} and therefore proves
Theorem~\ref{mainresult}. Note that because of the way
we constructed the annuli we have $|z_i-z_j|>64h$ for $i\neq j$
and therefore $B_{31h}(z_i)$ and $B_{31h}(z_j)$ are always a
distance of more than $h$ apart.
\end{proof}

\end{document}

%% file: figu00.pstex_t
\begin{picture}(0,0)%
\epsfig{file=figu00.pstex}%
\end{picture}%
\setlength{\unitlength}{3947sp}%
\begingroup\makeatletter\ifx\SetFigFont\undefined%
\gdef\SetFigFont#1#2#3#4#5{%
  \reset@font\fontsize{#1}{#2pt}%
  \fontfamily{#3}\fontseries{#4}\fontshape{#5}%
  \selectfont}%
\fi\endgroup%
\begin{picture}(3187,2342)(5903,-4457)
\put(7426,-3736){\makebox(0,0)[lb]{\smash{\SetFigFont{12}{14.4}{\familydefault}{\mddefault}{\updefault}\special{ps: gsave 0 0 0 setrgbcolor}2h\special{ps: grestore}}}}
\put(8101,-3886){\makebox(0,0)[lb]{\smash{\SetFigFont{12}{14.4}{\familydefault}{\mddefault}{\updefault}\special{ps: gsave 0 0 0 setrgbcolor}4h\special{ps: grestore}}}}
\end{picture}

%% file: fig000.pstex_t
\begin{picture}(0,0)%
\epsfig{file=fig000.pstex}%
\end{picture}%
\setlength{\unitlength}{3947sp}%
\begingroup\makeatletter\ifx\SetFigFont\undefined%
\gdef\SetFigFont#1#2#3#4#5{%
  \reset@font\fontsize{#1}{#2pt}%
  \fontfamily{#3}\fontseries{#4}\fontshape{#5}%
  \selectfont}%
\fi\endgroup%
\begin{picture}(3660,1887)(3553,-4036)
\put(5326,-3013){\makebox(0,0)[lb]{\smash{\SetFigFont{12}{14.4}{\familydefault}{\mddefault}{\updefault}\special{ps: gsave 0 0 0 setrgbcolor}y\special{ps: grestore}}}}
\put(5440,-3638){\makebox(0,0)[lb]{\smash{\SetFigFont{12}{14.4}{\familydefault}{\mddefault}{\updefault}\special{ps: gsave 0 0 0 setrgbcolor}Cat(y)\special{ps: grestore}}}}
\put(4474,-3979){\makebox(0,0)[lb]{\smash{\SetFigFont{12}{14.4}{\familydefault}{\mddefault}{\updefault}\special{ps: gsave 0 0 0 setrgbcolor}A rescaling of Cat(y)\special{ps: grestore}}}}
\put(4929,-2445){\makebox(0,0)[lb]{\smash{\SetFigFont{12}{14.4}{\familydefault}{\mddefault}{\updefault}\special{ps: gsave 0 0 0 setrgbcolor}$N_{\theta_0}(y)$\special{ps: grestore}}}}
\end{picture}

%% file: figu0.pstex_t
\begin{picture}(0,0)%
\epsfig{file=figu0.pstex}%
\end{picture}%
\setlength{\unitlength}{3947sp}%
\begingroup\makeatletter\ifx\SetFigFont\undefined%
\gdef\SetFigFont#1#2#3#4#5{%
  \reset@font\fontsize{#1}{#2pt}%
  \fontfamily{#3}\fontseries{#4}\fontshape{#5}%
  \selectfont}%
\fi\endgroup%
\begin{picture}(3092,1440)(8093,-2994)
\end{picture}

%% file: figu1.pstex_t
\begin{picture}(0,0)%
\epsfig{file=figu1.pstex}%
\end{picture}%
\setlength{\unitlength}{3947sp}%
\begingroup\makeatletter\ifx\SetFigFont\undefined%
\gdef\SetFigFont#1#2#3#4#5{%
  \reset@font\fontsize{#1}{#2pt}%
  \fontfamily{#3}\fontseries{#4}\fontshape{#5}%
  \selectfont}%
\fi\endgroup%
\begin{picture}(4878,2047)(9229,-4091)
\put(10006,-3114){\makebox(0,0)[lb]{\smash{\SetFigFont{9}{10.8}{\familydefault}{\mddefault}{\updefault}\special{ps: gsave 0 0 0 setrgbcolor}z\special{ps: grestore}}}}
\put(10164,-3744){\makebox(0,0)[lb]{\smash{\SetFigFont{9}{10.8}{\familydefault}{\mddefault}{\updefault}\special{ps: gsave 0 0 0 setrgbcolor}y'\special{ps: grestore}}}}
\put(9954,-2169){\makebox(0,0)[lb]{\smash{\SetFigFont{9}{10.8}{\familydefault}{\mddefault}{\updefault}\special{ps: gsave 0 0 0 setrgbcolor}$z^a$\special{ps: grestore}}}}
\put(10899,-3219){\makebox(0,0)[lb]{\smash{\SetFigFont{9}{10.8}{\familydefault}{\mddefault}{\updefault}\special{ps: gsave 0 0 0 setrgbcolor}$z^b$\special{ps: grestore}}}}
\put(13370,-3103){\makebox(0,0)[lb]{\smash{\SetFigFont{9}{10.8}{\familydefault}{\mddefault}{\updefault}\special{ps: gsave 0 0 0 setrgbcolor}z'\special{ps: grestore}}}}
\put(9954,-2746){\makebox(0,0)[lb]{\smash{\SetFigFont{9}{10.8}{\familydefault}{\mddefault}{\updefault}\special{ps: gsave 0 0 0 setrgbcolor}y\special{ps: grestore}}}}
\end{picture}

%% file: figu2.pstex_t
\begin{picture}(0,0)%
\epsfig{file=figu2.pstex}%
\end{picture}%
\setlength{\unitlength}{3947sp}%
\begingroup\makeatletter\ifx\SetFigFont\undefined%
\gdef\SetFigFont#1#2#3#4#5{%
  \reset@font\fontsize{#1}{#2pt}%
  \fontfamily{#3}\fontseries{#4}\fontshape{#5}%
  \selectfont}%
\fi\endgroup%
\begin{picture}(3112,2128)(5080,-9682)
\put(7864,-8047){\makebox(0,0)[lb]{\smash{\SetFigFont{8}{9.6}{\familydefault}{\mddefault}{\updefault}\special{ps: gsave 0 0 0 setrgbcolor}$z_2'$\special{ps: grestore}}}}
\put(7065,-9077){\makebox(0,0)[lb]{\smash{\SetFigFont{8}{9.6}{\familydefault}{\mddefault}{\updefault}\special{ps: gsave 0 0 0 setrgbcolor}$z_2$\special{ps: grestore}}}}
\put(5975,-8887){\makebox(0,0)[lb]{\smash{\SetFigFont{8}{9.6}{\familydefault}{\mddefault}{\updefault}\special{ps: gsave 0 0 0 setrgbcolor}$z_1$\special{ps: grestore}}}}
\put(5317,-7849){\makebox(0,0)[lb]{\smash{\SetFigFont{8}{9.6}{\familydefault}{\mddefault}{\updefault}\special{ps: gsave 0 0 0 setrgbcolor}$z_1'$\special{ps: grestore}}}}
\put(5701,-8236){\makebox(0,0)[lb]{\smash{\SetFigFont{12}{14.4}{\familydefault}{\mddefault}{\updefault}\special{ps: gsave 0 0 0 setrgbcolor}$\gamma_1$\special{ps: grestore}}}}
\put(7426,-8461){\makebox(0,0)[lb]{\smash{\SetFigFont{12}{14.4}{\familydefault}{\mddefault}{\updefault}\special{ps: gsave 0 0 0 setrgbcolor}$\gamma_2$\special{ps: grestore}}}}
\end{picture}

%% file: figu3.pstex_t
\begin{picture}(0,0)%
\epsfig{file=figu3.pstex}%
\end{picture}%
\setlength{\unitlength}{3947sp}%
\begingroup\makeatletter\ifx\SetFigFont\undefined%
\gdef\SetFigFont#1#2#3#4#5{%
  \reset@font\fontsize{#1}{#2pt}%
  \fontfamily{#3}\fontseries{#4}\fontshape{#5}%
  \selectfont}%
\fi\endgroup%
\begin{picture}(3905,2404)(7038,-8094)
\put(7501,-6511){\makebox(0,0)[lb]{\smash{\SetFigFont{12}{14.4}{\familydefault}{\mddefault}{\updefault}\special{ps: gsave 0 0 0 setrgbcolor}$\gamma_0$\special{ps: grestore}}}}
\put(9001,-6361){\makebox(0,0)[lb]{\smash{\SetFigFont{12}{14.4}{\familydefault}{\mddefault}{\updefault}\special{ps: gsave 0 0 0 setrgbcolor}$\gamma_1$\special{ps: grestore}}}}
\put(9526,-6886){\makebox(0,0)[lb]{\smash{\SetFigFont{12}{14.4}{\familydefault}{\mddefault}{\updefault}\special{ps: gsave 0 0 0 setrgbcolor}$\gamma_2$\special{ps: grestore}}}}
\put(10126,-7486){\makebox(0,0)[lb]{\smash{\SetFigFont{12}{14.4}{\familydefault}{\mddefault}{\updefault}\special{ps: gsave 0 0 0 setrgbcolor}$\gamma_3$\special{ps: grestore}}}}
\end{picture}